\documentclass[aps,preprint,showkeys,amsmath,amssymb,
aps,
]{revtex4-1}

\usepackage{graphicx}
\usepackage{dcolumn}
\usepackage{bm}
\usepackage{color}
\usepackage[utf8]{inputenc}
\usepackage[toc,page]{appendix}

\usepackage{mathtools}
\usepackage{empheq}
\usepackage[skins,theorems]{tcolorbox} 
\tcbset{highlight math style={enhanced,
		colframe=black,colback=white,arc=4pt,boxrule=1pt}}

\usepackage{algpseudocode}
\usepackage{algorithm}
\usepackage[english]{babel}

\usepackage{etoolbox}
\usepackage{amsmath}

\newcommand{\algstrut}[1][\algruledefaultfactor]{\vrule width 0pt
	depth .25\baselineskip height #1\baselineskip\relax}
\newcommand*{\algrule}[1][\algorithmicindent]{\hspace*{0.5em}\vrule\algstrut
	\hspace*{\dimexpr#1+1.0em}}

\makeatletter
\newcount\ALG@printindent@tempcnta
\def\ALG@printindent{%
	\ifnum \theALG@nested>0
	\ifx\ALG@text\ALG@x@notext
	\else
	\unskip
	\ALG@printindent@tempcnta=1
	\loop
	\algrule[\csname ALG@ind@\the\ALG@printindent@tempcnta\endcsname]%
	\advance \ALG@printindent@tempcnta 1
	\ifnum \ALG@printindent@tempcnta<\numexpr\theALG@nested+1\relax
	\repeat
	\fi
	\fi
}%

\patchcmd{\ALG@doentity}{\noindent\hskip\ALG@tlm}{\ALG@printindent}{}{\errmessage{failed to patch}}

\AtBeginEnvironment{algorithmic}{\lineskip0pt}

\begin{document}

\date{\today}

\keywords{chaos control, transient chaos, time series.}


\title{Controlling two-dimensional chaotic transients with the safety function}

\author{Rub{\'e}n Cape{\'a}ns}
\affiliation{Nonlinear Dynamics, Chaos and Complex Systems Group, Departamento de  F\'isica, Universidad Rey Juan Carlos, M\'ostoles, Madrid, Tulip\'an s/n, 28933, Spain}

\author{Miguel A.~F. Sanju{\'a}n}
\email{Corresponding author: miguel.sanjuan@urjc.es}
\affiliation{Nonlinear Dynamics, Chaos and Complex Systems Group, Departamento de  F\'isica, Universidad Rey Juan Carlos, M\'ostoles, Madrid, Tulip\'an s/n, 28933, Spain}

\begin{abstract}
	
	In this work we deal with the Hénon and the Lozi map for a choice of parameters where they show transient chaos. Orbits close to the chaotic saddle behave chaotically for a while to eventually escape to an external attractor. Traditionally, to prevent such an escape,  the partial control technique has been applied. This method stands out for considering disturbances (noise) affecting the map and for finding a special region of the phase space, called \textit{the safe set}, where the control required to sustain the orbits is small. However, in this work we will apply a new approach of the partial control method that has been recently developed. This new approach is based on finding a special function called \textit{the safety function} which allows to automatically find the minimum control necessary to avoid the escape of the orbits. Furthermore, we will show the strong connection between the safety function and the classical safe set. To illustrate that, we will compute for the first time, safety functions for the two-dimensional Hénon and Lozi maps, where we also show the strong dependence of this function with the magnitude of disturbances affecting the map, and how this change drastically impacts the controlled orbits. 

\end{abstract}

\maketitle

\section{Introduction}

Chaotic transient behavior \cite{Ottbook,TransientChaos} usually arises when  due to the change of a parameter of the system, a chaotic attractor collides with his own boundary basin causing a boundary crisis. The chaotic set becomes  a nonattracting chaotic set (i.e chaotic repellor or a chaotic saddle),and almost all orbits in the neighborhood of the nonattracting chaotic set are free to escape to an external attractor. If in addition we consider that the dynamics are affected for some disturbances (noise), all trajectories sooner or later eventually escape.

Traditionally, to avoid the escape of the orbits and sustain them around the nonattracting chaotic set, we have used the classical partial control method. This technique is applied on maps and is based on finding certain special set in the phase space, called the safe set, where the control needed is small. Unlike other control methods ( \cite{Biological}, Schwartz and Triandaf \cite{Schwartz}, and Dhamala and Lai \cite{Dhamala}), partial control is designed from the ground up to deal with disturbances. Furthermore, it takes into account the magnitude of the disturbances, to find the best safe set. One of the more remarkable result is that the control used to sustain the orbits in the safe set, is always smaller than the  magnitude of disturbances affecting the orbits.

However, in this work we will not use the classical approach of the partial control method based on finding a safe set. Instead, we will use the new approach based on the computation of a special functions called the safety function from which we can obtain the minimum control necessary to sustain the orbits and also the minimum safe set. In this sense, this new approach is a generalization of the classical partial control method.

In the next sections we will illustrate the application of this new approach based on the computation of the safety function. First with the one-dimensional tent map, and then, for the first time, with the two-dimensional Hénon and Lozi map.

\section{The partial control: from the classical approach to the new approach}

The first step to apply the partial control technique is to define a region $Q$ in the phase space containing the chaotic saddle. In this region $Q$, we assume that the dynamics can be described by the following map:

\begin{equation}
\begin{array}{l}
q_{n+1}=f(q_n)+\xi_n+u_n, \hspace{1cm}\text{with} \hspace{0.4cm}
|\xi_n| \leq \xi_0, \hspace{0.3cm} |u_n| \leq u_0 < \xi_0,\\
\nonumber
\end{array}
\end{equation}
where $q$ describes the state vector of the system, $\xi_n$ is the  disturbance affecting the map and we assume it is limited by an upper bound value $\xi_0$.  The control $u_n$ is applied every iteration of the map with the knowledge of $f(q_n) + \xi_n$. This control is also limited by the upper bound value $u_0$, which can be previously fixed by the controller. 

The  points $q\in Q$ that can be sustained in $Q$ without exceed the control $u_0$, defines the safe set (which is a subset of $Q$). To do that, there is an algorithm called the Sculpting Algorithm \cite{Automatic}. This algorithm  takes as input, the map $f$ defined in $Q$ and the values $\xi_0$ and $u_0$, and compute the corresponding safe set. However if the value $u_0$ introduced is too small, no safe set exists. Therefore to find the minimum $u_0$ value, we have to gradually increasing  $u_0$  until the safe set appears. This task can be rather unpractical since $u_0$ is manually set for every computation. In a recent work ~\cite{SafetyFun}, it was presented a new approach of partial control, that automatically computes the minimum $u_0$ value. This approach is based on the computation of a function in the region $Q$ called the safety function $U$. The algorithm for the computation is shown in the Appendix. This algorithm takes as input the map $f$ defined in $Q$ and the value $\xi_0$, and computes the corresponding safety function. The minimum of this function corresponds to the minimum $u_0$ possible. The minimum safe set is the set of points $q\in Q$ that satisfy $U(q)=u_0$.

To show an example of how the safety function works, we will use the well known tent map slope-three for a choice of the parameter where transient chaos is present.

\begin{equation}
\label{}
x_{n+1} = \left\{
\begin{array}{ll}
\mu x_n + \xi_n +u_n      & \mathrm{\;for\ } x_n \le \frac{1}{2} \\
\mu(1-x_n) + \xi_n +u_n    & \mathrm{\;for\ } x_n > \frac{1}{2} \\
\end{array}
\right.
\end{equation}

\vspace{1cm}

\begin{figure}
	\includegraphics [trim=0cm 0cm 0cm 0cm, clip=true,width=0.9\textwidth]{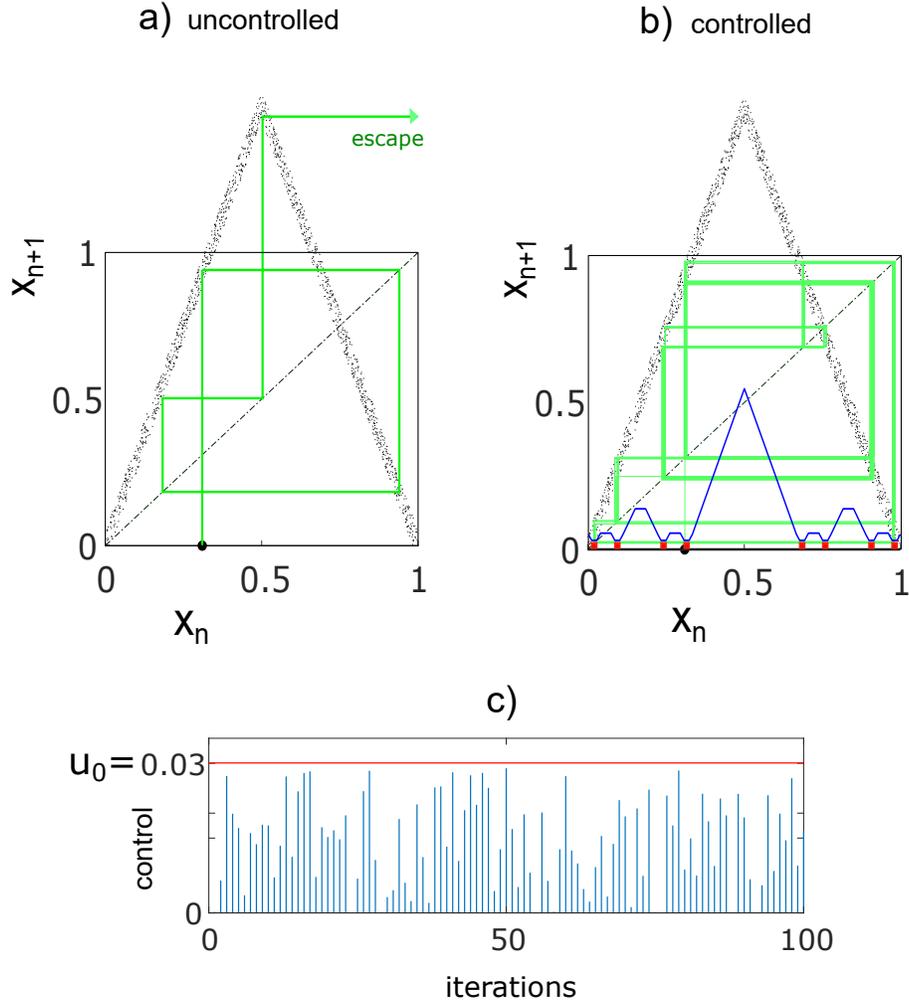}\\
	\centering
	\caption{\textbf{Partial control method.} The slope-three tent map is represented in this figure. The map is affected by a uniform disturbance bounded by $\xi_0=0.05$. The small dots help to visualize the intensity and distribution of the disturbance. (a) An uncontrolled orbit that escapes from the interval $Q=[0,1]$ after a few iterations. (b) In blue, the safety function that has 8 minima with value $u_0=0.03$. These 8 minima defines the safe set, which is represented with the red pieces at the bottom. In green, a controlled orbit is shown. This controlled orbit starts in the point $x_0=0.3$ that belongs to the safe set. At every iteration of the map, the orbit is forced to pass through the safe set to remain forever in $Q=[0,1]$. (c) Controls $|u_n|\leq u_0=0.03$ applied during the first $100$ iterations of the map.}
	\label{1}
\end{figure}

This map with $\mu=3$, exhibits transient chaos in the interval $Q=[0,1]$ (see Fig.~\ref{1}a). We consider that orbits of this map are affected by disturbances $\xi_n \leq \xi_0=0.05$. then, we compute the corresponding safety function shown in Fig.~\ref{1}b in blue line. As shown,the safety function has $8$ minima with the value $u_0=0.03$. The minimum safe set corresponds with the location of this minima indicated in Fig.~\ref{1}b by the small red pieces. Finally, orbits starting in the the safe set can be sustained inside it by applying every iteration a control $|u_n| \leq u_0=0.03$. In Fig.~\ref{1}b, 100 iterations of a controlled orbit (green line) is represented and the corresponding 100  $|u_n|$ controls in Fig.~\ref{1}c.

\section{Partial control applied to the Hénon map and the Lozi map}

In previous works \cite{Automatic, Ecology, Lorenz, Parametric}, safe sets have been computed from 1D, 2D and 3D maps, by applying the classical algorithm described in \cite{Asymptotic,Automatic}. However,  computing the safe sets by means of the safety function  has been mainly shown with one-dimensional maps as the tent map presented before. 

In this work, we want to show the computation of the safety functions in two-dimensional maps. The extra dimension just adds more computation since now the safety function will be a two-dimensional surface $U(q_x,q_y)$. However the algorithm to compute the safety function (see the Appendix) remains the same since it applies to any dimension. 

To illustrate the control method, we choose the well known Hénon map and the Lozi map for a choice of parameters where they show transient chaos. In both cases we will follow the same procedure: First define the region $Q$ containing the  chaotic saddle and set the upper bound of disturbance $xi_0$ affecting the map in $Q$. Next, compute the corresponding safety function. Finally obtain the minimum safe set and control the orbits to remain in it.

\subsection{Application to the H\'enon map}

In 1976 the French astronomer Michel Hénon introduced the map later named after him,  defined as:
\begin{equation}
\begin{array}{l}
x_{n+1}=a-b y_n- x_n^2  \\
y_{n+1}=x_n.   \\
\end{array}
\end{equation}

Hénon proved that this is the most general form of quadratic maps, which  shows transient chaos for a wide range of parameters $a$ and $b$. Here we have have chosen the parameter values $a=6$ and $b=0.4$. For these values, the trajectories with initial conditions in the square $Q=[-4,4] \times [-4,4]$ have a short chaotic transient, before finally escaping this region towards infinity. 

The corresponding Hénon map including the disturbance and the control is:

\begin{equation}
\begin{array}{ l }
x_{n+1}=a-b y_n-x_n^2+\xi^x_n+u^x_n  \\
y_{n+1}=x_n +\xi^y_n+u^y_n. \\ 
\end{array}
\end{equation}

where we choose the bound of disturbance $\xi_0=0.20$ so that $\parallel\xi^x_n, \xi^y_n\parallel \leq \xi_0$. The control is also limited so that $\parallel u^x_n, u^y_n\parallel \leq u_0$, where the minimum $u_0$ possible, will be defined by the minimum value of the safety function. Below this value, no safe set exists.

\begin{figure}
	\includegraphics [trim=0cm 0cm 0cm 0cm, clip=true,width=1\textwidth]{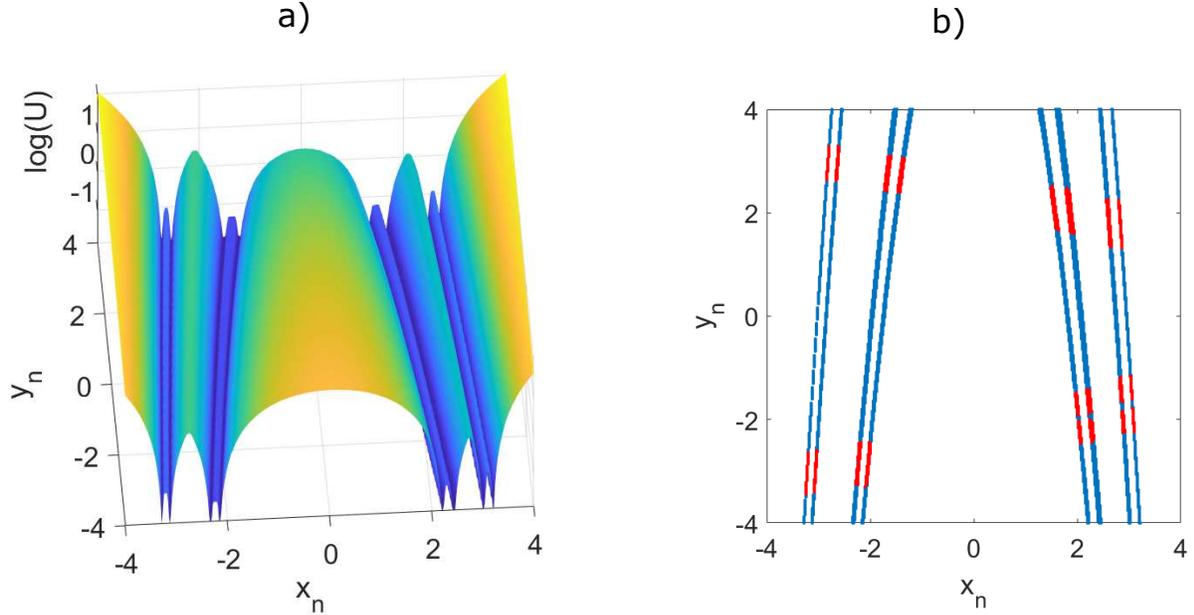}\\
	\centering
	\caption{\textbf{The 2D safety function for the H\'{e}non map.} The bound of a disturbance is $\xi_0=0.20$. a) The safety function was computed in the square $Q=[-4,4] \times [-4,4]$ using a grid of $1000\times1000$ points. The safety function takes 18 iterations to converge and it has the minimum value $\min(U)=0.15=u_0$. The function is represented logarithmic to enhance the visualization. b) The points $q\in Q$ that satisfy $U(q)=0.15$ define the safe set (in blue).  We also represent 10000 iterations of a controlled orbit marked with the red dots.}
	\label{2}
\end{figure}

The corresponding safety function is shown in Fig.~\ref{2} a, where it has been plotted in logarithmic scale  for a better visualization. For this case the minimum value is found to be $\min(U)=0.15=u_{0}$. The points $q\in Q$ with $U(q)=0.15$ defines the minimum safe set, which is represented in the  Fig.~\ref{2} b. 

Finally we use the safe set to control the orbits inside it. Every iteration of the map, the control $(u^x_n,u^y_n)$ with $\parallel u^x_n, u^y_n\parallel \leq 0.15$ is applied to put the orbit back in the nearest safe point. A controlled orbit of 10000 iterations is represented in Fig.~\ref{2} b by the red dots.  

\subsection{Application to the Lozi map}

A piecewise linear version of the Hénon map was given by the French mathematician René Lozi in the form:

\begin{equation}
\begin{array}{l}
x_{n+1}=1-a|x_n|+b y_n  \\
y_{n+1}=x_n.   \\
\end{array}
\end{equation}

For the choice of parameters $a=2$, $b=0.5$ the orbits in the square $Q=[-4,4] \times [-4,4]$ behave chaotic for a while to eventually escape from it. In order to avoid this escape we will applied the partial control technique following the same steps as with the Hénon map. In this case the controlled maps is the following:

\begin{equation}
\begin{array}{l}
x_{n+1}=1-a|x_n|+b y_n +\xi^x_n+u^x_n \\
y_{n+1}=x_n +\xi^y_n+u^y_n   \\\\
\end{array}
\end{equation}

with $\parallel\xi^x_n, \xi^y_n\parallel \leq \xi_0$ and $\parallel u^x_n, u^y_n\parallel \leq u_0$.

\begin{figure}
	\includegraphics [trim=0cm 0cm 0cm 0cm, clip=true,width=1\textwidth]{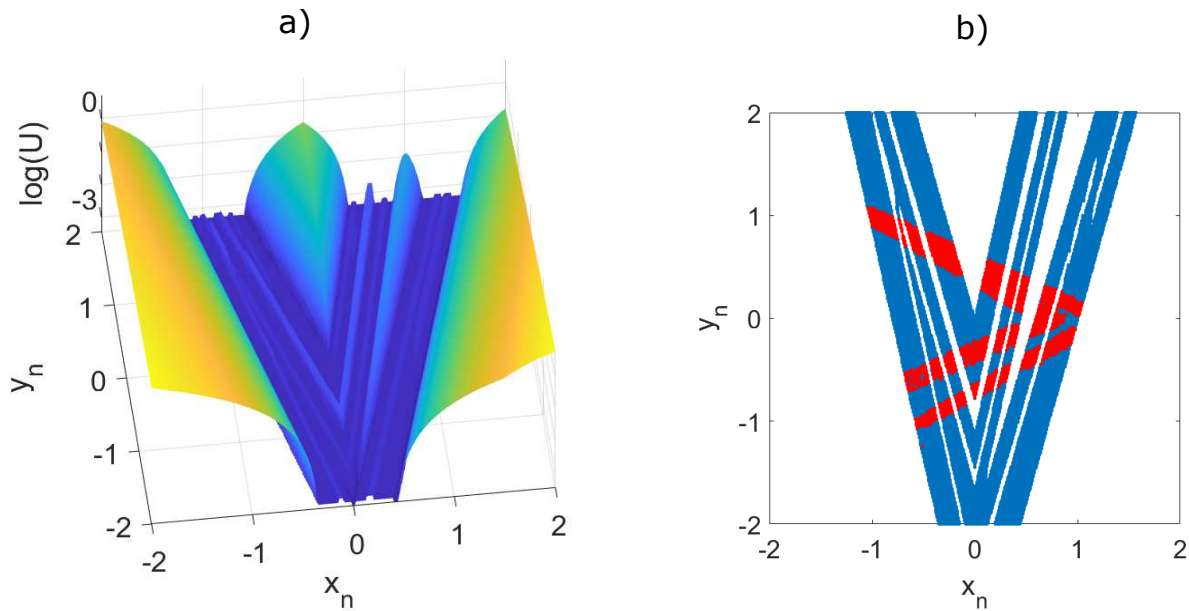}\\
	\centering
	\caption{\textbf{The 2D safety function for the Lozi map.} The bound of a disturbance selected is $\xi_0=0.050$. a) The safety function was computed in the square $Q=[-4,4] \times [-4,4]$ using a grid of $1000\times1000$ points. The safety function takes 23 iterations to converge and it has the minimum value $\min(U)=0.035=u_0$. The function is represented logarithmic to improve the visualization. b) By taking the points $q\in Q$ that satisfy $U(q)=0.035$ we obtain the safe set (blue). We also represent 10000 iterations of a controlled orbit marked with the red dots.}
	\label{3}
\end{figure}

To compute an example, we consider the upper bound of disturbance $\xi_0=0.050$. Then we compute the safety function $U$ shown in Fig.~\ref{3} a, where the function is plotted logarithmic for a better visualization. The minimum of $U$ is found to be $\min(U)=0.035=u_0$. Therefore, the minimum safe set is the set of point $q\in Q$ for which $U(q)=0.035$. This safe set is shown in the in Fig.~\ref{3} b, where we also draw 10000 iterations of a controlled orbit (red points).

\section{Safe set variation with \boldmath$\xi_0$}

\begin{figure}
	\includegraphics [trim=3cm 0cm 2cm 0cm, clip=true,width=1.0\textwidth]{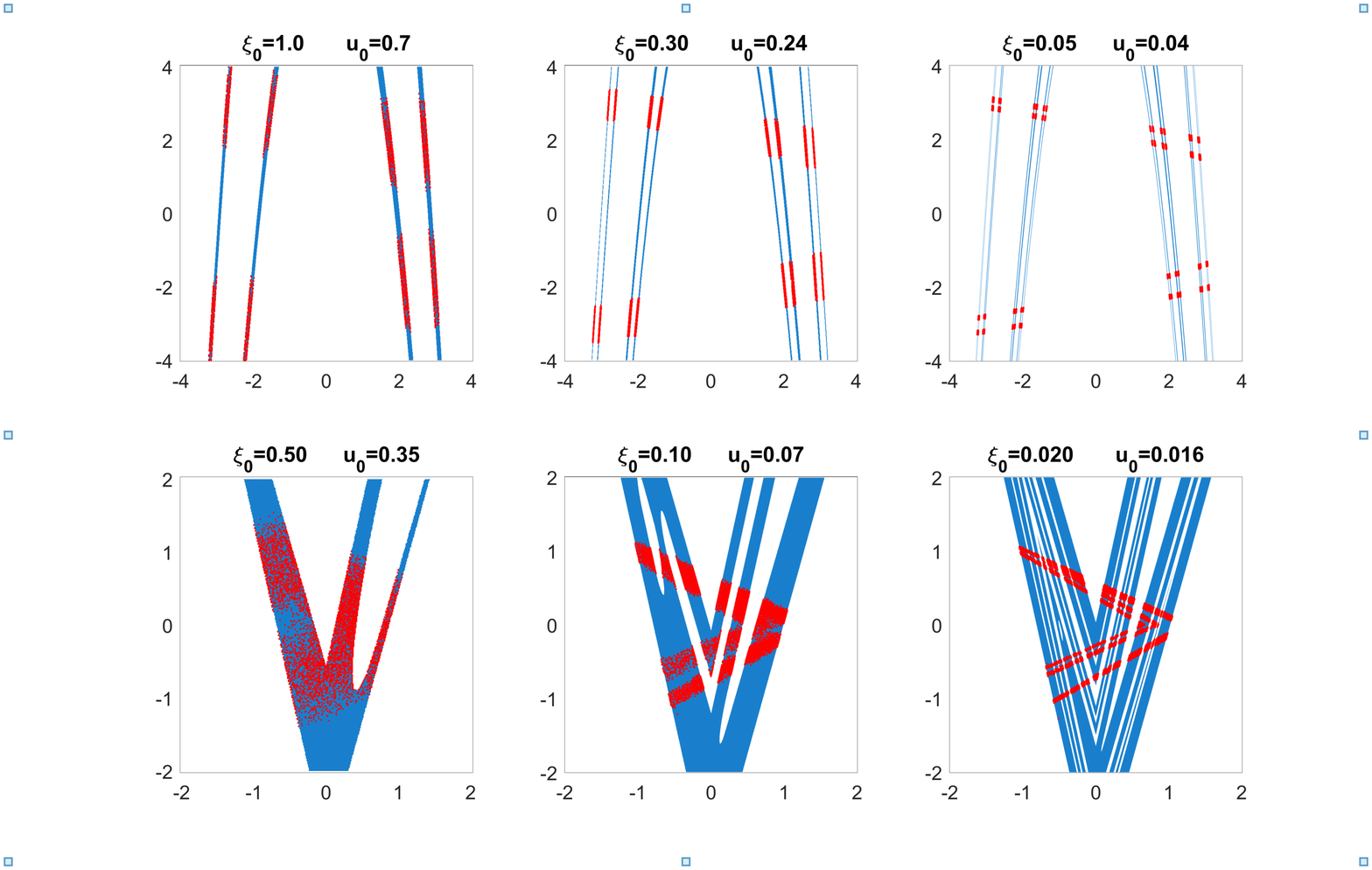}\\
	\centering
	\caption{\textbf{Safe sets for different \boldmath$\xi_0$.} The top figure correspond to safe sets computed for the Hénon map affected by different upper disturbance values $\xi_0$, that decreases from left to right. On the bottom, the  safe sets corresponding to the Lozi map. All the safe sets are minimum with the minimum bound of control $u_0 < \xi_0$  indicated on the top. Below this control value, no safe set exists. Note that the safe sets computed with biguer $\xi_0$, consist of a few and fat strips, while the safe sets computed with smaller $\xi_0$, are made of many and thin strips.  This is very clear in the case of the Hénon safe sets where we have 4,8 and 16 strips respectively. The red dots plotted over each safe set correspond to 10000 iterations of a controlled orbit.}
	\label{4}
\end{figure}

The partial control method takes advantage of the fractal structure of the non-attracting chaotic set responsible for the chaotic transients. As shown in the two-dimensional maps used here, the safe sets resembles the coarse grained structure of the stable manifold of the chaotic saddle. The grain size is mainly determined by the magnitude of disturbance affecting the map. Smaller values of $\xi_0$ leads to finer safe sets like a Cantor set. In Fig  ~\ref{4}, different safe sets for the Hénon and the Lozi map has been computed.
Top safe sets correspond to the Hénon map and bottom figures to the Lozi map.
In both cases the bound of disturbance $\xi_0$ affecting the map, decreases from left to right figures.

It each safe set of Fig  ~\ref{4}, it has been also plotted a controlled orbits consisting of 10000 iterations. We can observe  that any starting orbit  in the safe set, quickly converges to a smaller region called the asymptotic safe set \cite{Asymptotic}, that resembles the coarse grained structure of  the chaotic saddle.  Once the controlled orbit enters in the asymptotic safe set, it remains inside forever..  

Finally, we want to point out that, although the safe sets computed here were the smaller possible, there is no restriction to compute safe sets with bigger values of $u_0$. These safe sets will be a fattened version of the minimum safe set.  This feature can be of interest if we want that the controlled orbit visit more points of the region $Q$, at the expense of applying bigger controls  $|u_n|\leq u_0$.

\section{Conclusions}

In this work we have shown the application of the new approach of partial control technique based on the safety functions. This technique is applied to  maps showing transient chaos with the goal to avoid the escape of the orbits from the non-attractive chaotic set. For the first time, we applied this technique to the two-dimensional Hénon map, and the Lozi map, both affected by disturbance that is considered bounded. In each case, we define a region $Q$ containing the chaotic saddle, where we compute the safety function. Then we have extracted the minimum safe set, where the orbits can be sustaining using a minimum control bound. Finally we have shown how the minimum safe sets change depending on the bound of disturbance affecting the map, which has a drastic impact in the controlled orbits.

\begin{algorithm}[H]
	\renewcommand{\thealgorithm}{}
	\floatname{algorithm}{}
	\caption{\textbf{Appendix: The safety function algorithm}}\label{euclid}
	
	\begin{algorithmic}[0]
		\Statex
		\Statex 
		\textbf{Notation}:\\
		 	$i\equiv$ index of the starting point $q[i],~ i=1:N$ where $N=$ total number of grid points in $Q$.\\
		$s\equiv$ index of the disturbance $\xi[s],~ s=1:M$ where $M=$ number of disturbed images.\\
		$j\equiv$ index of the arrival point $q[j],~ j=1:N$.\\
		Map with the index notation: $~q_{n+1}=f(q_n) + \xi_n +u_n$ $~~~\rightarrow~~~$  $q[j]=\,f\big(q[i]\big) +\xi[s]+ u\,[i,s,j]$.\\
		Computation of the safety function $U$:

			\Statex
			\State- Initially set $\;U_0\,[j]=0, \; \forall j=1:N, \; \; \; \; k=0.$
			\State
			\hspace{4cm}\While {$U_{k+1}\neq U_k$}
			\State
			\For{$i=1$ to  $N\;$}
			\For{$s=1$ to  $M_i\;$}
			\For{$j=1$ to  $N\;$}
			\State $\;u\,[i,s,j]= \Big|\,f\big(q[i]\big)+\xi[s]- \;q\,[j]\,\Big|$
			\Comment{\parbox[t]{.33\linewidth}{\scriptsize{Distance between the disturbed image $f\big(q[i]\big)+\xi[s]$ and the arrival point $q\,[j]$. Note that the $u\,[i,s,j]$ values remain unchanged every iteration of the while loop so compute them once and save them.}}}
			\State
			\State $\;u^*\,[i,s,j]=\max\limits_{j}\,\big(\,u[i,s,j],\, U_k[j]\,\big)$
			\EndFor
			\State $\;u^{**}\,[i,s]=\min\limits_{j}\,\big(\,u^*\,[i,s,j]\,\big)$
			\EndFor
			\State $U_{k+1}[i]=\max\limits_{s} \big(\,u^{**}\,[i,s]\,\big)$
			\EndFor
			\State  $k=k+1$\normalsize
			\State
			\EndWhile
			\Statex
		
		\State Compact formula: $U_{k+1}[i]=\max\limits_{1\leq s\leq M_i}\Big(\min\limits_{1\leq j\leq N}\big(\max\limits_{j}\, (\,u[i,s,j],\, U_k[j]\,)\,\big)\,\Big)$ \\

		\Statex \\
	\end{algorithmic}
\end{algorithm}


\begin{acknowledgments}
This work was supported by the Spanish State Research Agency (AEI) and the European Regional Development Fund (FEDER) under Project No. PID2019-105554GB-100.
\end{acknowledgments}

\section*{Compliance with ethical standards}
\textbf{Conflict of interest} The authors declare that they have no conflict of interest concerning the publication of this manuscript.

\section*{Data availability}
All data generated or analysed during this study have been obtained through numerical simulations that are included in this article

\end{document}